  % This is the Plain TeX file for the paper
  % Non-Hausdorff groupoids,  by Ruy Exel

  \magnification 1200

%\input /Users/exel/Casa/z/p/nhausdorf/WEIRD/head.tex

  % FONTS

  \newcount\fontset
  \fontset=1
  \def \dualfont#1#2#3{\font#1=\ifnum\fontset=1 #2\else#3\fi}

  \dualfont\bbfive{bbm5}{cmbx5}
  \dualfont\bbseven{bbm7}{cmbx7}
  \dualfont\bbten{bbm10}{cmbx10}

  \font \eightbf = cmbx8
  \font \eighti = cmmi8 \skewchar \eighti = '177
  \font \eightit = cmti8
  \font \eightrm = cmr8
  \font \eightsl = cmsl8
  \font \eightsy = cmsy8 \skewchar \eightsy = '60
  \font \eighttt = cmtt8 \hyphenchar\eighttt = -1

  \font \sixi = cmmi6 \skewchar \sixi = '177
  \font \sixrm = cmr6
  \font \sixsy = cmsy6 \skewchar \sixsy = '60
  \font \tensc = cmcsc10

  \scriptfont \bffam = \bbseven
  \scriptscriptfont \bffam = \bbfive
  \textfont \bffam = \bbten

  \newskip \ttglue

  \def \eightpoint {\def \rm {\fam0 \eightrm }%
  \textfont0 = \eightrm
  \scriptfont0 = \sixrm \scriptscriptfont0 = \fiverm
  \textfont1 = \eighti
  \scriptfont1 = \sixi \scriptscriptfont1 = \fivei
  \textfont2 = \eightsy
  \scriptfont2 = \sixsy \scriptscriptfont2 = \fivesy
  \textfont3 = \tenex
  \scriptfont3 = \tenex \scriptscriptfont3 = \tenex
  \def \it {\fam \itfam \eightit }%
  \textfont \itfam = \eightit
  \def \sl {\fam \slfam \eightsl }%
  \textfont \slfam = \eightsl
  \def \bf {\fam \bffam \eightbf }%
  \textfont \bffam = \bbseven
  \scriptfont \bffam = \bbfive
  \scriptscriptfont \bffam = \bbfive
  \def \tt {\fam \ttfam \eighttt }%
  \textfont \ttfam = \eighttt
  \tt \ttglue = .5em plus.25em minus.15em
  \normalbaselineskip = 9pt
  \def \MF {{\manual opqr}\-{\manual stuq}}%
  \let \sc = \sixrm
  \let \big = \eightbig
  \setbox \strutbox = \hbox {\vrule height7pt depth2pt width0pt}%
  \normalbaselines \rm }

  % REFERENCE CONTROL

  % Define \showlabel, \showref and/or \showcit if you wanto to see labels,
  % references, or citations, eg.
  % \def \showlabel{} \def \showref{} \def \showcit{}

  \newcount \secno \secno = 0
  \newcount \stno \stno = 0
  \newcount \eqcntr \eqcntr= 0

  \def \ifn #1{\expandafter \ifx \csname #1\endcsname \relax }

  \def \track #1#2#3{\ifn{#1}\else {\tt\ [#2 \string #3] }\fi}

  \def \advseqnumbering {\global \advance \stno by 1 \global \eqcntr=0}

  \def \current {\number \secno \ifnum \number \stno = 0 \else
    .\number \stno \fi }

  \def \laberr#1#2{\message{*** RELABEL CHECKED FALSE for #1 ***}
      RELABEL CHECKED FALSE FOR #1, EXITING.
      \end}

  \def \syslabel#1#2{%
    \ifn {#1}%
      \global \expandafter 
      \edef \csname #1\endcsname {#2}%
    \else
      \edef\aux{\expandafter\csname #1\endcsname}%
      \edef\bux{#2}%
      \ifx \aux \bux \else \laberr{#1=(\aux)=(\bux)} \fi
      \fi
    \track{showlabel}{*}{#1}}

  \def \subeqmark #1 {\global \advance\eqcntr by 1
    \edef\aux{\current.\number\eqcntr}
    \eqno {(\aux)}
    \syslabel{#1}{\aux}}

  \def \eqmark #1 {\advseqnumbering
    \eqno {(\current)}\syslabel{#1}{\current}}

  \def \label #1 {\syslabel{#1}{\current}}

  \def \lcite #1{(#1\track{showcit}{$\bullet$}{#1})}

  \def \cite #1{[{\bf #1}\track{showref}{\#}{#1}]}

  \def \scite #1#2{{\rm [\bf #1\track{showref}{\#}{#1}{\rm \hskip 0.7pt:\hskip 2pt #2}\rm]}}

  % HEADER

 \def \Headlines #1#2{\nopagenumbers
    \advance \voffset by 2\baselineskip
    \advance \vsize by -\voffset
    \headline {\ifnum \pageno = 1 \hfil
    \else \ifodd \pageno \tensc \hfil \lcase {#1} \hfil \folio
    \else \tensc \folio \hfil \lcase {#2} \hfil
    \fi \fi }}

  \def \Date #1 {\footnote {}{\eightit Date: #1.}}

  % CONTROL SEQUENCES

  \def \lcase #1{\edef \auxvar {\lowercase {#1}}\auxvar }
  \def \section #1{\global\def \SectionName{#1}\stno = 0 \global
\advance \secno by 1 \bigskip \bigskip \goodbreak \noindent {\bf
\number \secno .\enspace #1.}\medskip \noindent \ignorespaces}

  \long \def \sysstate #1#2#3{%
    \advseqnumbering
    \medbreak \noindent 
    {\bf \current.\enspace #1.\enspace }{#2#3\vskip 0pt}\medbreak }
  \def \state #1 #2\par {\sysstate {#1}{\sl }{#2}}
  \def \definition #1\par {\sysstate {Definition}{\rm }{#1}}
  \def \remark #1\par {\sysstate {Remark}{\rm }{#1}}

  % Examples
  % \sysstate ...{Theorem}.{font}{Text}
  % \state .......Theorem.........Text\par
  % \definition ..................Text\par

  \def \proof {\medbreak \noindent {\it Proof.\enspace }}
  \def \proofend {\ifmmode \eqno \square \else \hfill \square
\looseness = -1 \medbreak \fi }

  \def \$#1{#1 $$$$ #1}
  \def \=#1{\buildrel \hbox{\sixrm #1} \over =}

  \def \Item #1{\smallskip \item {{\rm #1}}}
  \newcount \zitemno \zitemno = 0

  \def \izitem {\zitemno = 0}
  \def \zitemplus {\global \advance \zitemno by 1 \relax}
  \def \rzitem{\romannumeral \zitemno}
  \def \rzitemplus {\zitemplus \rzitem}
  \def \zitem {\Item {{\rm(\rzitemplus)}}}
  \def \zitemmark #1 {\syslabel{#1}{\rzitem}}

  \newcount \nitemno \nitemno = 0
  
  \def \nitem {\global \advance \nitemno by 1 \Item {{\rm(\number\nitemno)}}}

  \newcount \aitemno \aitemno = -1
  \def \boxlet#1{\hbox to 6.5pt{\hfill #1\hfill}}
  \def \iaitem {\aitemno = -1}
  \def \aitemconv{\ifcase \aitemno a\or b\or c\or d\or e\or f\or g\or
h\or i\or j\or k\or l\or m\or n\or o\or p\or q\or r\or s\or t\or u\or
v\or w\or x\or y\or z\else zzz\fi}
  \def \aitem {\global \advance \aitemno by 1\Item {(\boxlet \aitemconv)}}
  \def \aitemmark #1 {\syslabel{#1}{\aitemconv}}

  \newcount \footno \footno = 1
  \newcount \halffootno \footno = 1
  \def \footcntr {\global \advance \footno by 1
  \halffootno =\footno
  \divide \halffootno by 2
  $^{\number\halffootno}$}
  \def \fn#1{\footnote{\footcntr}{\eightpoint#1\par}}

  % STANDARD DEFINITIONS

  \def \C {{\bf C}}
  \def \<{\left \langle \vrule width 0pt depth 0pt height 8pt }
  \def \>{\right \rangle }  
  
  \def \and {\hbox {,\quad and \quad }}

  \def \imply {\mathrel{\Rightarrow}}
  \def \for #1{,\quad \forall\,#1}
  \def \square {\hbox {$\sqcap \!\!\!\!\sqcup $}}
  
  \def \stress #1{{\it #1}\/}
  \def \inv {^{-1}}
  \def \*{\otimes}
  \def\chi{{supp'}}

  \newcount \bibno \bibno = 0
  \def \newbib #1{\global\advance\bibno by 1 \edef #1{\number\bibno}}

  \def \bibitem #1#2#3#4{\smallskip \item {[#1]} #2, ``#3'', #4.}

  \def \references {
    \begingroup
    \bigskip \bigskip \goodbreak
    \eightpoint
    \centerline {\tensc References}
    \nobreak \medskip \frenchspacing }

  %%%%%%%%%%%%%%%%%%%%%%%%%
  \input pictex

  \newcount\ax
  \newcount\ay
  \newcount\bx
  \newcount\by
  \newcount\dx
  \newcount\dy
  \newcount\vecNorm
  \newcount\pouquinho \pouquinho = 200

  \def\beginmypicture{ \begingroup \noindent \hfill \beginpicture}
  \def\endmypicture{\endpicture \hfill\null \endgroup}

  \def\myarrow#1#2#3#4{\arrow <0.15cm> [0.25,0.75] from #1 #2 to #3 #4 }%
  \def\morph#1#2#3#4{%
    \ax = #1
    \ay = #2
    \bx = #3
    \by = #4
    \dx = \bx \advance \dx by -\ax
    \dy = \by \advance \dy by -\ay
    \vecNorm = \dx 
    \ifnum\vecNorm<0 \vecNorm=-\vecNorm \fi
    \advance \vecNorm by \ifnum\dy>0 \dy \else -\dy \fi
    \multiply \dx by \pouquinho \divide \dx by \vecNorm
    \multiply \dy by \pouquinho \divide \dy by \vecNorm
    \advance \ax by \dx
    \advance \bx by -\dx
    \advance \ay by \dy
    \advance \by by -\dy
    \myarrow{\number\ax}{\number\ay}{\number\bx}{\number\by}}
    %%%%%%%%%%%%%%%%%%%%%%%%%%%%%%%

  \def\Gz{G^{(0)}}
  \def\CRG{C^*_r(G)}
  \def\CZGZ{C_0\big(\Gz\big)}
  \def\germ[#1,#2,#3]{[\kern1pt #2,#3]}
  \def\ng#1#2{[\kern1pt #1,#2\kern1pt]}
  \def\R{{\bf R}}
  
  \def\halfbox#1{\hbox to 0.45\hsize{$\scriptstyle \bullet$ #1\hfill}}
  \def\dbitem#1#2{\medskip\noindent\hskip 10pt\hbox{\halfbox{#1}\halfbox{#2}}}

  \newbib\RenaultClaire
  \newbib\Connes
  \newbib\CK
  \newbib\Exel
  \newbib\infinoa
  \newbib\topfree
  \newbib\Kumjian
  \newbib\KPR
  \newbib\PatBook
  \newbib\RenaultBook
  \newbib\Renault
  \newbib\RS

  \def\titletextOne{NON-HAUSDORFF GROUPOIDS}

  \Headlines {\titletextOne} {R.~Exel}

  \null\vskip -1cm
  \centerline{\bf \titletextOne} 
  \footnote{\null}
  {\eightrm 2000 \eightsl Mathematics Subject Classification:
  \eightrm 
  Primary   46L55, % Noncommutative dynamical systems
  secondary 22A22. % Topological groupoids.
  }

  \centerline{\tensc 
    R.~Exel\footnote{*}{\eightpoint Partially supported by
CNPq.}}\footnote{\null}
  {\eightrm Keywords: non-Hausdorff groupoids,  essentially principal
groupoids.}
  \Date{21 Nov 2009}

  \midinsert 
  \narrower \narrower
  \eightpoint \noindent We present examples of non-Hausdorff,
\'etale, essentially principal groupoids for which three results, known
to hold in the Hausdorff case, fail.  These results are: 
  \quad (A)
the subalgebra of continuous functions on the unit space is maximal
abelian within the reduced groupoid C*-algebra, 
  (B) every nonzero ideal of the reduced groupoid C*-algebra has a
nonzero intersection with the subalgebra of continuous functions on
the unit space, and
  (C) the open support of a normalizer
  is a bissection.

  \endinsert

\section{Introduction} 
  This paper is concerned with \stress{\'etale groupoids}
  \cite{\RenaultBook,\RenaultClaire,\PatBook,\Renault,\Exel}.
  A topological groupoid $G$ is said to be \'etale
if its unit space $\Gz$ is locally compact and Hausdorff, and the
range map ``$r$" (and consequently also the source map ``$s$") is a
local homeomorphism.

One may or may not assume the global topology of $G$ to be Hausdorff
but, while non-Hausdorff \stress{topological spaces} may be safely
ignored in numerous applications of Topology, non-Hausdorff groupoids
do occur in many essential situations, such as the holonomy groupoid
of a foliation \cite{\Connes} or the groupoid of germs of a
pseudogroup of local homeomorphisms on a topological space
\scite{\Renault}{Section 3}.

Therefore, rather than dismissing non-Hausdorff groupoids as a
nuisance, it is highly desirable to embrace them in the general
theory.

  An \'etale groupoid $G$ is said to be \stress{principal} if its
\stress{isotropy group bundle}, namely
  $$
  G' := \{\gamma\in G: s(\gamma) = r(\gamma)\},
  $$
  coincides with the unit space $\Gz$, and it is said to be
\stress{essentially principal} if the \stress{interior} of $G'$
coincides with $\Gz$.  Principal groupoids correspond to 
\stress{free} group actions while the essentially principal ones
correspond to \stress{topologically free} actions,  hence the
relevance of these concepts.
  
Among the important consequences of the property of being essentially
principal, in the Hausdorff case,  the following stand out: 
  \bigskip\item{(A)} $\CZGZ$ is maximal abelian within the reduced
groupoid C*-algebra $\CRG$ \scite{\Renault}{4.2}.
  \medskip\item{(B)} Every nonzero ideal of $\CRG$ has a nonzero
intersection with $\CZGZ$ (see the appendix for a precise statement).
  \medskip\item{(C)} If $u$ is a normalizer of $\CZGZ$ within $\CRG$,
then the \stress{open support} of $u$, namely
  $$
  \chi(u) = \{\gamma\in G: u(\gamma)\neq 0\},
  $$
  is a bissection \scite {\Renault}
{Proposition 4.7}.

\bigskip These results underlie mainstream developments in the
theory of C*-algebras:  (B) is related to
  uniqueness theorems for Cuntz-Krieger algebras
  \scite{\CK}{2.15}, \scite{\infinoa}{13.2}
  and graph algebras \cite{\KPR}, 
  as well as to results on reduced crossed products by partial group actions
  \scite{\topfree}{2.6}, 
  while (A) and (C) are related to Cartan subalgebras \cite{\Kumjian, \Renault}.

It should be stressed that (A),  (B)  and (C) are only  known to hold under
the assumption that $G$ is Hausdorff!

In trying to embrace non-Hausdorff groupoids within the general
theory, I (and quite likely many other people) have spent a lot of
energy in the effort to generalize the above facts beyond the
Hausdorff case.  After having failed to do so I have found examples of
non-Hausdorff \'etale groupoids which provide counter-examples for all
of the above statements.  In what follows we shall discuss these
examples in detail.

Our first example is related to an example by G. Skandalis \cite{\RS} built
with a different purpose, namely of exhibiting a minimal foliation
whose C*-algebra is not simple.

I would like to thank Jean Renault for many fruitful discussions, and
for bringing Skandalis' example to my attention.
I would also like to thank Alcides Buss for many interesting
discussions while I was searching for the second example below.

\section{The first example}
Consider the following subsets of $\R^2$:
  $$
  \matrix{
  X  & = &  [-1, 1]\times \{0\}, \cr\cr
  Y & = & \{0\}\times [-1, 1],\cr\cr
  Z & = & X \cup Y.\hfill}
  $$
  Clearly $Z$ is invariant under the action of the subgroup
$H\subseteq GL_2({\bf R})$ generated by
  $$
  \sigma_x = \pmatrix{ -1 & 0 \cr 0 & 1}
  \and
  \sigma_y = \pmatrix{ 1 & 0 \cr 0 & -1}.
  $$

  Let $G$ be the groupoid of germs for the action of $H$ on $Z$ (see
Section (3) of \cite{\Renault} for the definition of the groupoid of
germs for a given pseudogroup).
  As is the case for every groupoid of germs, $G$ is essentially
principal \scite{\Renault}{3.4}.

We shall adopt a slightly simplified notation in relation to
\cite{\Renault}, namely the germ of the transformation $\varphi$ at the
point $x$ will be denoted by $\ng\varphi x$, as opposed to Renault's
notation $[y,\varphi,x]$, where $y=\varphi(x)$.

In the present case it is interesting to observe that,
  $$
  \matrix{
  \ng{\sigma_y}{x} = \ng{I}{x}, & \quad
  \ng{\sigma_x\sigma_y}{x} = \ng{\sigma_x}{x},\cr\cr
  \ng{\sigma_x}{y} = \ng{I}{y}, & \quad
  \ng{\sigma_x\sigma_y}{y} = \ng{\sigma_y}{y},}
  $$
  for all $x\in X^* := X{\setminus}\{0\}$, and all $y\in Y^* :=
Y{\setminus}\{0\}$,
  where ``$I$\kern1pt" stands for the identity map, and we denote the
zero vector of $\R^2$ simply by ``$0$".
  We therefore see that $G$ consists of the following distinct
elements:

  \dbitem
  {$\ng{I}{x}$, for $x\in X^*$,}
  {$\ng{I}{y}$, for $y\in Y^*$,}
  \dbitem
  {$\ng{\sigma_x}{x}$, for $x\in X^*$,}
  {$\ng{\sigma_y}{y}$, for $y\in Y^*$,}
  \dbitem
  {$\ng{I}{0}$,}
  {$\ng{\sigma_x}{0}$,}
  \dbitem
  {$\ng{\sigma_y}{0}$,}
  {$\ng{\sigma_x\sigma_y}{0}$.}

\bigskip Observe that the isotropy group bundle $G'$ is formed by the
last three elements listed above, in addition to the units.

Recall that a \stress{bisection} is a subset of $G$ restricted to
which both the range and source maps are injective.  Consider the
following open bisections of $G$:

  \dbitem
  {$U_1 = \big\{ \ng{I}{z}: z\in Z\big\} = \Gz$,}
  {$U_x = \big\{ \ng{\sigma_x}{z}: z\in Z\big\}$,}
  \dbitem
  {$U_y = \big\{ \ng{\sigma_y}{z}: z\in Z\big\}$,}
  {$U_{xy} = \big\{ \ng{\sigma_x\sigma_y}{z}: z\in Z\big\}$.}

\bigskip\noindent 
  Let $f_1,f_x,f_y,f_{xy}\in C_c(G)$ (for the definition of $C_c(G)$
see \cite{\Connes}, \cite{\PatBook} or \scite{\Exel}{3.9}) be the
characteristic function of $U_1,U_x,U_y$, and $U_{xy}$, respectively.
Finally put
  $$
  f = f_1 - f_x - f_y + f_{xy}.
  $$
By direct computation one checks that 
  \dbitem
  {$f(\ng{I}{0}) = 1$,}
  {$f(\ng{\sigma_x}{0}) = -1$,}
  \dbitem
  {$f(\ng{\sigma_y}{0}) = -1 $,}
  {$f(\ng{\sigma_x\sigma_y}{0}) = 1$,}
  \medskip  \noindent and that $f$ vanishes on all other points of $G$.
  In particular notice that the support of $f$ (set of points where
$f$ does not vanish, no closure) is the set
  $$ 
  \big\{\ng{I}{0},
  \ng{\sigma_x}{0},
  \ng{\sigma_y}{0},
  \ng{\sigma_x\sigma_y}{0}\big\} = r\inv(\{0\}) = s\inv(\{0\}),
  \eqmark SinvZero
  $$
  which is contained in $G'$.

\state Proposition \label MainComputation For every $g\in C_c(G)$ one
has that
  $$
  g*f = f*g = \lambda(g)f,
  %  \Big(g(\ng I0)-g(\ng{\sigma_x}0)-g(\ng{\sigma_y}0)+
  %  g(\ng{\sigma_x\sigma_y}0)\Big)f
  $$
  where $\lambda(g)$ is the scalar given by $\lambda(g) = g(\ng
I0)-g(\ng{\sigma_x}0)-g(\ng{\sigma_y}0)+
    g(\ng{\sigma_x\sigma_y}0)$.

\proof
  Recall that for every $\gamma\in G$ one has
  $$
  (f*g)(\gamma) = \sum_{\alpha\beta=\gamma}f(\alpha)g(\beta)
  \for\gamma\in G.
  $$
  If the above sum is nonzero, then there exists at least one pair
$(\alpha,\beta)$ such that $\alpha\beta=\gamma$, and $f(\alpha)\neq0$.
As seen in \lcite{\SinvZero}, this implies that $r(\alpha)=0$, and
hence necessarily $r(\gamma) =0$, as well.  Therefore $f*g$ is
supported in $r\inv(\{0\})$.  A similar reasoning and the same
conclusion applies to $g*f$.

We leave it to the reader to compute $(f*g)(\gamma)$ and
$(g*f)(\gamma)$ for the four elements $\gamma$ in $r\inv(\{0\})$, after
what the result will become apparent.
  \proofend

The first conclusion to be drawn from the above result is:

\state Proposition Even though $G$ is essentially principal (every
groupoid of germs is essentially principal by \scite{\Renault}{3.4}),
there is a nonzero ideal $J\subseteq \CRG$ for which $J\cap \CZGZ =
\{0\}$.

\proof By \lcite{\MainComputation} one has that $J:=\C f$ is an ideal
in $\CRG$.  Since $f$ is not in $\CZGZ$, the intersection of
$J$ with $\CZGZ$ is trivial.  \proofend

The second conclusion is:

\state Proposition Even though $G$ is essentially principal,
one has that $\CZGZ$ is not maximal abelian within $\CRG$.

\proof
  If is enough to notice that by \lcite{\MainComputation} one has that
$f$ is a central element of $\CRG$, and hence commutes with every
element of $\CZGZ$, but $f$ is not in $\CZGZ$.
  \proofend

Since the support of $f$ is contained in $G'$, and in view of
\scite{\Renault}{4.2}, it is not surprising that $f$ commutes with
every element of $\CZGZ$.  
  % However something else seems to be taking place here since $f$ is
  % actually central.

  \font\sc = cmcsc9
  \def\g{\gamma}
  \def\ep{\varepsilon}
  \def\V{{\cal V}}
  \def\c{\subseteq}
  \def\rep{\pi}
  \def\Ker{{\rm Ker}}
  \def\CCG{C_c(G)}
  \def\a{\alpha}
  \def\vbar{\overline{v}}

\def\B{C_0\big(G^{(0)}\big)} \def\B{C_0(X)}
\def\A{C^*_r(G)}
\def\s{\sigma}
\section{Strange normalizers} Let $n$ be a positive integer and let
  $
  I=I_n=\{1, 2, \ldots, n\}
  $
  be seen as a discrete topological space.  On the product space
  $[0, 1]\times I$, we consider the equivalence relation ``$\sim$" according to which
  $$
  (0, i) \sim (0, j) \for i, j\in I,
  $$
  and such that no other pairs of points are related except for each point with itself.
The quotient topological
space
  $$
  X = \big([0, 1]\times I\big)/\sim
  $$
  therefore looks like a star with $n$ edges.  Incidentally, in the
special and very relevant case $n=4$, notice that $X$ is homeomorphic
to the space $Z$ of the previous section.

The equivalence class
of $(0,i)$ will be denoted simply by $0$, and if $t>0$, the
equivalence class of $(t,i)$, namely the singleton $\{(t,i)\}$, will
be denoted by $(t,i)$, by abuse of language.

Let $S_n$ be the group of permutations of $I$ and consider the action
of $S_n$ on $[0, 1]\times I$, where each $\s\in S_n$ acts by
  $$
  (t, i) \mapsto (t, \s(i)).
  $$ 
  The equivalence relation ``$\sim$" above is clearly left invariant
by this action so we get an action of $S_n$ on $X$.
  Considering the subgroup $A_n\subseteq S_n$ formed by all even
permutations, we may restrict the above  action to $A_n$, and we shall
let
  $$
  G=G(X, A_n)
  $$
  be the corresponding groupoid of germs.  The unit space of $G$ is
therefore homeomorphic to $X$ and we shall tacitly identify these from
now on.

The main technical result of this section is in order:

\state Theorem \label WeirdNormalizer Assuming that $n\geq4$, and
given any $\tau\in S_n$ (not necessarily in $A_n$), there exists a
unitary element $u$ in $\A$ such that for all $t>0$, and all $i,j\in
I$, one has 
  $$
  u(t,i,j) = 
  \left\{\matrix{ 1, & \hbox{ if } \tau(i) = j, \cr 0, & 
    \hbox{ otherwise.}
    \vrule height 13pt width 0pt
    }\right.
  $$
  Moreover 
  $
  u^*fu = f\circ \tau,
  $
  for every $f\in \B$.

\proof 
  % If $x\in X$ and $\s\in S_n$, recall that the germ of $\s$ at $x$
  % is denoted $[\s, x]$.
  Notice that if $x=(t,i)$, with $t>0$, and if $\s, \s'\in S_n$, then
  $$
  [\s, (t, i)] =   [\s', (t, i)] 
  \iff \s(i) = \s'(i).
  $$
  In other words, the germ of $\s$ at $(t,i)$ depends only on
$j:=\s(i)$.
  We may therefore denote this germs simply by $(t, i, j)$.
  On the other hand, it is easy to see that
  $$
  [\s, 0] =   [\s', 0] 
  \iff \s = \s'.
  $$
  We may then describe $G$ as being the set
  $$
  \big\{(t,i,j): t\in(0,1],\ i,j\in I \big\} \quad \cup \quad
\big\{[\s,0]: \s\in A_n\big\}.
  $$

For each $\s\in A_n$, let $U_\s=\{[\s, x]: x\in X\}$ be the canonical
bissection associated to $\s$.  Since $U_\s$ is compact one has that
its characteristic function, here denoted $1_\s$, is an
element of $\A$ which is easily seen to be unitary.  Moreover, the
correspondence
  $$
  \s\in A_n \mapsto 1_\s\in \A
  $$
  is a unitary representation of $A_n$ in $\A$, which therefore
integrates to a *-homomorphism 
  $$
  \phi: C^*(A_n) \to \A.
  $$
  Given a generic element
  $$
  a = \sum_{\s\in A_n}a_\s \delta_\s\in C^*(A_n),
  \subeqmark FormOfElement  % \eqno{(\dagger)}
  $$
  and any $(t,i,j)\in G$, with $t> 0$, observe that
  $$
  \phi(a) (t,i,j) =
  \sum_{\s\in A_n}a_\s 1_\s(t,i,j) =
  \sum_{\s(i)=j}a_\s.
  \subeqmark ExplicitCalculation  % \eqno{(\ddagger)}
  $$

  Changing subjects slightly, consider the representation $\pi$ of
$S_n$ on the Hilbert space $\C^n$, where each $\s\in S_n$ is mapped
to the unitary operator $\pi(\s)$ defined on the cannonical basis
$\{e_i\}_{i\in I}$ of $\C^n$ by
  $$
  \pi(\s)e_i = e_{\s(i)}.
  $$
  Denote by $\tilde\pi$ the corresponding integrated representation of
$C^*(S_n)$ on $\C^n$, and observe that for each $a$ as in \lcite{\FormOfElement},
one has that
  $$
  \<\tilde\pi(a)e_i, e_j\> =
  \sum_{\s\in A_n}a_\s \<\pi(\s)e_i,e_j\> =
  \sum_{\s(i)=j}a_\s,
  $$
  so we have by \lcite{\ExplicitCalculation}  that
  $$
  \phi(a) (t,i,j) =  \<\tilde\pi(a)e_i, e_j\>,
  \subeqmark StarEqtn % \eqno{(\star)}
  $$
  for all $a\in C^*(A_n)$, all $t\in (0, 1]$, and all $i,j\in I$.

Assuming that $n\geq4$, one may prove that the commutant in ${\cal
B}(\C^n)$ of both $\tilde\pi\big(C^*(S_n)\big)$ and
$\tilde\pi\big(C^*(A_n)\big)$ coincide with the set of all matrices of
the form
  $$
  \pmatrix{
  z & y & y & \ldots & y \cr
  y & z & y & \ldots & y \cr
  y & y & z & \ldots & y \cr
  \vdots & \vdots & \vdots & \ddots & \vdots \cr
  y & y & y & \ldots & z \cr
  },
  $$
  where $z,y\in\C$.  The crucial point in doing this is that $A_n$
acts \stress{bi-transitively} on $I$, meaning that given
$i_1\neq i_2$ and $j_1\neq j_2$, there exists $\s\in A_n$ such that
$\s(i_1)=j_1$, and $\s(i_2)=j_2$.  Incidentally this is not true for $n<4$.

By the double commutant Theorem we conclude that 
$\tilde\pi\big(C^*(S_n)\big) = \tilde\pi\big(C^*(A_n)\big)$.
Given $\tau\in S_n$,  as in the statement,  we therefore have that $\pi(\tau)\in
\tilde\pi(C^*(A_n))$, so there exists some $v\in C^*(A_n)$ such that
$\tilde\pi(v) = \pi(\tau)$.  Absent any K-theoretic obstructions we
may assume that $v$ is unitary.

The element $u$ of which the statement speaks is 
the unitary
element $u:=\phi(v)\in \A$.  To see that it  satisfies the required
conditions notice that 
  for all $t>0$, and $i,j\in I$, we have that
  $$
  u(t,i,j) = 
  \phi(v)(t,i,j) \={\lcite{\StarEqtn}}
  \<\tilde\pi(v)e_i, e_j\> =
  \<\pi(\tau)e_i, e_j\> =
  \<e_{\tau(i)}, e_j\>,
  % \left\{\matrix{ 1, & \hbox{ if } \tau(i) = j, \cr 0, & \hbox{ otherwise.}}\right.
  $$
  proving the first assertion.
  It follows that the \stress{open support} of $u$, namely
  $$
  \chi(u) = \{\gamma\in G: u(\gamma)\neq 0\},
  \subeqmark DefineOpenSup
  $$
  (cf. \cite {\Renault}) consists precisely of all germs
$(t,i,\tau(i))$, where $t> 0$, besides a few other germs at $0$.
Therefore the range of any $\gamma\in \chi(u)$ coincides whith the
image of its source under the action of $\tau$.  From this it
immediately follows that
  $$
  u^*fu = f\circ \tau
  \for f\in \B.
  \proofend
  $$

The relevance of this result is in relation to \scite {\Renault}
{Proposition 4.7}, where it is proved that the open support of a
normalizer of $\B$ in $\A$ is a bissection.  In the present
non-Hausdorff situation this fails:

\state Proposition Let $n\geq4$, let $\tau\in S_n\setminus A_n$, and
let $u$ be given as in \lcite{\WeirdNormalizer}. Then $u$ is a
normalizer of $\CZGZ$ within $\CRG$ but,  even though
$G$ is essentially principal, the open
support of $u$ is not a bissection.

\proof
  Recall from the proof of \lcite{\WeirdNormalizer} that $u=\phi(v)$,
where $v$ is a unitary element in $C^*(A_n)$ such that
$\tilde\pi(v)=\pi(\tau)$.
  Assuming by contradiction that $\chi(u)$ is a bissection, and
noticing that the germs $[\s,0]$ all have range and source equal to
$0$, we deduce that there is at most one $\s\in A_n$ for which 
$u([\s,0])\neq0$.

  Write $v = \sum_{\s\in A_n} a_\s \delta_\s$, as in \lcite{\FormOfElement},
so that 
$u = \sum_{\s\in A_n} a_\s 1_\s$, and hence 
  $$
  u([\s,0]) = a_\s \for \s\in A_n.
  $$
  We then see that there is only one $\s\in A_n$ for which $a_\s\neq0$,
which implies that $v=a_\s \delta_\s$.  Consequently
  $$
  \pi(\tau) = \tilde\pi(v) = a_\s \pi(\s),
  $$
  which contradicts the fact that $\tau\notin A_n$.  \proofend

Another dilemma presented by this example is related to the program
initiated by Kumjian in \cite{\Kumjian}, and recently continued by
Renault in \cite{\Renault}, attempting to classify Cartan sub-algebras
of C*-algebras.  Given a commutative subalgebra $B$ of a C*-algebra
$A$ and a normalizer $u\in N(B)$, Kumjian \cite{\Kumjian} showed the
existence of a partial homeomorphism $\theta_u$ of the spectrum $X$ of
$B$, such that for all $x\in X$, with $(u^*u)(x)\neq0$, one has
  $$
  (u^*bu)(x) = (u^*u)(x) b\big(\theta_u(x)\big)
  \for b\in B.
  $$
  When $A$ is the reduced C*-algebra of an essentially principal,
Hausdorff, \'etale groupoid $G$, and $B=C_0\big(\Gz\big)$, Renault
showed \cite{\Renault} that one can reconstruct $G$ from the inclusion
``$B\c A$", as the germs of the partial homeomorphisms $\theta_u$,
where $u$ ranges in the set of all normalizers.

In the present example, if one attempted to reconstruct $G$ from the
inclusion
  ``$C_0(X)\c \CRG$"
  using the above method, the presence of the strange normalizers $u$
above would lead us to consider the germ at zero of every $\tau\in
S_n$ by the last assertion of \lcite{\WeirdNormalizer}, but
the isotropy group $G(0)$ is only as big as $A_n$!

\section{Appendix (\it the intersection property for ideals in essentially
principal,  Hausdorff groupoids)}
  In this section we prove result (B) stated in the introduction.
Although this result has been used in several contexts under various
guises (see the introduction for some references), it seems not to
have appeared in the literature in quite the general form we have in
mind.  We begin with some elementary considerations about
representations of commutative C*-algebras.

  Let $X$ be a locally compact Hausdorff space and let $\rep$ be a
representation of $C_0(X)$ on a Hilbert space $H$.  As any ideal of
$C_0(X)$, the kernel of $\rep$ must be of the form $C_0(U)$, for some
open set $U\c X$.

\definition Given a representation $\rep$ of $C_0(X)$, with
$\Ker(\rep)=C_0(U)$, we will refer to
$X{\setminus} U$ as the \stress{support of\/ $\rep$}.

\state Lemma
  \label EasyCommutLemma
  If $\rep$ is a representation of $C_0(X)$ on a Hilbert space $H$,
and if  $x$ lies in the support of $\rep$, then $|f(x)|\leq \|\rep(f)\|$,
for every $f$ in $C_0(X)$.

\proof Left to the reader.  \proofend

Given a groupoid $G$, for every $x\in\Gz$ we denote by $G(x)$
  ({\it cf.} \scite{\RenaultBook}{I.1.1})
  the \stress{isotropy group at $x$}, namely
  $$
  G(x) = \{\g\in G: s(\g)=r(\g) = x\}.
  $$
  Obviously $x\in G(x)$, but in case $G(x)=\{x\}$ we say that $x$ \stress{has
no isotropy}.

The following result gives the key inequality from which the next Theorem
will be deduced.  It is roughly based on \scite{\RenaultBook}{II.4.4}.

\state Lemma
  \label KeyInequality  Let $G$ be an \'etale, Hausdorff groupoid and 
  let $\rep$ be a representation of $C^*(G)$ on a Hilbert space $H$.
Suppose in addition that we are given $x\in \Gz$ such that
  \izitem
  \zitem $x$ has no isotropy,
  % that is,  $r(\g) = s(\g)=x \imply \g=x$,
  \zitem $x$ lies in the support of\/ $\rep|_{C_0(\Gz)}$.
  \medskip\noindent
Then for every $f\in\CCG$, one has that
  $
  |f(x)|\leq \|\rep(f)\|.
  $

\proof Let $\V$ be the collection of all open neighborhoods of $x$
within $\Gz$.  We will view $\V$ as a directed set under the order
relation
  $$
  v\leq w \iff v\supseteq w
  \for v,w\in \V.
  $$
  For each $v\in \V$, choose $v_1,v_2\in\V$, relatively compact, and such
that $\vbar_2\c v_1\c \vbar_1\c v$.  By Uryshon's Lemma let
  $$
  g_v:\Gz\to[0,1]
  $$
  be a continuous function whose restriction to $\vbar_2$ is
identically equal to $1$, and which vanishes off $v_1$.  The support
of $g_v$ is contained in $\vbar_1$, which is compact, so
$g_v\in C_c(\Gz)$.
  % Moreover, by (ii) and \lcite{\EasyCommutLemma} one has that
  % $$
  % 1 = g_v(x) \leq \|\rep(g_v)\|\leq\|g_v\| = 1,
  % $$
  % so $\|\rep(g_v)\|= 1$.

  We claim that there exists $\xi_v\in H$,
with $\|\xi_v\|=1$, and such that 
  $$
  \rep(g_v)\xi_v=\xi_v.
  $$
  In order to prove it, use Uryshon's Lemma again to produce a
continuous function $h_v:\Gz\to[0,1]$, vanishing off $v_2$ and such
that $h_v(x)\neq0$.
  Observe that, since $g_v$ is identically equal to $1$ on $v_2$, we
have that
  $$
  g_vh_v=h_v.
  \subeqmark gvhvhv
  $$
  From (ii) and \lcite{\EasyCommutLemma} it follows that 
  $$
  0<|h_v(x)| \leq \|\rep(h_v)\|,
  $$
  so $\rep(h_v)\neq0$, and one may pick
$\eta_v\in H$ such that $\|\rep(h_v)\eta_v\|=1$.  Setting
$\xi_v=\rep(h_v)\eta_v$, we have that
  $$
  \rep(g_v)\xi_v = 
  \rep(g_v)\rep(h_v)\eta_v =
  \rep(g_vh_v)\eta_v \={(\gvhvhv)}
  \rep(h_v)\eta_v = 
  \xi_v, 
  $$
  proving our claim.
  We next claim  that
  $$
  \lim_{v\in \V}
  \<\vrule height9pt width0pt
  \rep(f)\xi_v,\xi_v\> = f(x)
  \for f\in \CCG.
  \subeqmark InnProfForIneq
  $$

Without loss of generality we will suppose that there are $K,U\c G$,
such that $K$ is compact, $U$ is an open bisection, $K\c U$, and $f$ vanishes
outside of $K$. We will further denote by $\a_U$ the
homeomorphism from $s(U)$ to $r(U)$ given by
  $$
  \a_U\big(s(\g)\big) = r(\g)
  \for \g\in U.
  $$
  The proof of \lcite{\InnProfForIneq} will be broken up in the
following three cases:
  \iaitem
  \aitem $x\notin s(U)$, 
  \aitem $x\in s(U)$, and $\a_U(x)\neq x$,
  \aitem $x\in s(U)$, and $\a_U(x)=x$.

\def\proofunder#1{\bigskip\noindent {\sc Proof of
  \lcite{\InnProfForIneq} under} (#1):}

\proofunder a
 Noticing that $x\notin s(K)$, there exists some $v_0\in\V$, with
$v_0\cap s(K)=\emptyset$.  For every $v\c v_0$, one then has that
  $$
  (fg_v)(\g) =f(\g)g_v(s(\g)) =0
  \for\g \in G, 
  $$
  because either $\g\notin K$, or $s(\g)\in s(K)$.  Therefore
  $\rep(f)\xi_v = \rep(fg_v)\xi_v = 0$, proving that the left-hand
side of \lcite{\InnProfForIneq} vanishes.  Observing that $x\notin K$
(or else $x=s(x)\in s(K)\c s(U)$), we see that the 
right-hand side of \lcite{\InnProfForIneq} also vanishes.

\proofunder b
  Let $A$ and $B$ be pairwise disjoint open subsets of $\Gz$ such that
$x\in A$ and $\a_U(x)\in B$.  Setting $v_0 = A\cap \a_U\inv(B)$,
notice that $x\in v_0$, and that $v_0\cap \a_U(v_0) = \emptyset$.  For
every $v\c v_0$ we have that
  $$
  (g_v^*fg_v)(\g) =
  \overline{g_v\big(r(\g)\big)} \; f(\g)\; g_v\big(s(\g)\big) 
  % = 0 
  \for\g \in G.
  $$
  If the above is nonzero for some $\g$, then $\g\in U$ and both
$s(\g)$ and $r(\g)$ lie in $v$.  Therefore
  $$
  r(\g) = \alpha_U\big(s(\g)\big)\in v\cap\alpha_U(v) \c
  v_0\cap\alpha_U(v_0) = \emptyset,
  $$
  which is impossible.  So $g_v^*fg_v=0$, and hence 
  $$ 
  \<\rep(f)\xi_v,\xi_v\> =
  \<\rep(f)\rep(g_v)\xi_v,\rep(g_v)\xi_v\> =
  \<\rep(g_v^*fg_v)\xi_v,\xi_v\> = 0,
  $$
  again proving the left-hand side of \lcite{\InnProfForIneq} to
vanish.  As for the right-hand side notice that $x\notin U$,
because otherwise 
  $$
  \alpha_U(x) =   \alpha_U\big(s(x)\big) = r(x) = x,
  $$
  which is in not in accordance with (b).  Thus $f(x)=0$, and
\lcite{\InnProfForIneq} is verified in the present case.

\proofunder c
  Given that $\a_U(x)=x$, there exists $\g\in U$ such that
$s(\g)=r(\g)=x$.  Thus
  $$
  \g\in G(x) =
  %  (isotropy?) =
  \{x\},
  $$
  so $\g=x$, by (i), and hence $x\in U$.  As $G$ is assumed to be
Hausdorff, we have that $f$ is 
  continuous\fn{On a non-Hausdorff groupoid the accepted definition of
$\CCG$ (see \cite{\Connes}, \cite{\PatBook} or \scite{\Exel}{3.9})
includes functions which are discontinuous, so we truly need to assume
$G$ to be Hausdorff here.}
  so, given $\ep>0$, we may choose a neighborhood $v_0$ of $x$,
contained in $U\cap \Gz$, and such that
  $$
  y\in v_0\imply |f(x)-f(y)|\leq\ep.
  $$

For every $v\c v_0$ we have that
  $$
  (fg_v)(\g) =
  f(\g)\; g_v\big(s(\g)\big)
  \for\g \in G.
  $$
  If the above is nonzero for some $\g$, then $\g\in U$ and $s(\g)\in
v$.  Thus, both $\g$ and $s(\g)$ lie in $U$, and since these have the
same source, we deduce that $\g=s(\g)$, and hence that $\gamma\in v$.
It follows that $fg_v$ vanishes outside $v$ and, in particular,
$fg_v\in\CZGZ$.  On the other hand, for every $y\in v$, one has that
  $$
  \big|(fg_v)(y)-f(x)g_v(y)\big| = 
  \big|f(y)-f(x)\big|\big| g_v(y)\big| \leq
  \ep \|g_v\|=\ep,
  $$ 
  which gives $\|fg_v-f(x)g_v\|\leq\ep$.
  Therefore, for $v$ as above,
  $$
  \big|\<\rep(f)\xi_v,\xi_v\> - f(x)\big| =
  \big|\<\rep(f)\xi_v,\xi_v\> - \<f(x)\xi_v, \xi_v\>\big| \$=
  \big|\<\rep(fg_v)\xi_v,\xi_v\> - \<f(x)\rep(g_v)\xi_v, \xi_v\>\big| =
  \big|\<\rep\big(fg_v-f(x)g_v\big)\xi_v,\xi_v\> \big| \$\leq
  \|fg_v-f(x)g_v\|\|\xi_v\|\|\xi_v\| \leq\ep, 
  $$
  proving \lcite{\InnProfForIneq} under the last case.
We then finally get
  $$
  |f(x)|  \={(\InnProfForIneq)}
  \lim_{v\in \V}|\<\rep(f)\xi_v,\xi_v\>| \leq
  \|\rep(f)\|.
  \proofend
  $$

We are now ready to prove the precise form of result (B) stated in the
introduction.

\state Theorem
  Let $G$ be an \'etale, Hausdorff, essentially principal, second
countable groupoid.
  % such that the set of points $x\in\Gz$ with trivial isotropy is
  % dense\fn{Such as an essentially principal groupoid whose unit
  % space has the Baire property \scite{\Renault}{3.1}.}.
  \iaitem
  \aitem If $\rep$ is a representation of $\CRG$ such that $\rep$ is
faithful on $\CZGZ$, then $\rep$ is faithful.
  \aitem If $J$ is a nonzero ideal in $\CRG$, then $J\cap\CZGZ$ is nonzero.

\proof
  We address (a) first.  Since $\rep$ is assumed to be faithful on
$\CZGZ$, the support of $\rep|_{C_0(\Gz)}$ is the whole of $\Gz$.
Given $f\in\CCG$ one then has by \lcite{\KeyInequality} that
  $$
  |f(x)| \leq \|\rep(f)\|,
  \subeqmark ApplicationOfKey
  $$
  for every $x$ in $\Gz$ without entropy.  Employing
\scite{\Renault}{3.1} we see that the set of such
$x$'s is dense in $\Gz$, and since the restriction of $f$
to $\Gz$ is continuous\fn{Again this would not be guaranteed should we not have
assumed that $G$ is Hausdorff.}, we conclude that in fact 
\lcite{\ApplicationOfKey} holds for every $x\in\Gz$, so
  $$
  \sup_{x\in\Gz}|f(x)|\leq\|\rep(f)\|.
  \subeqmark ApplicationOfKeyWithSup
  $$

  Let $E$ be the standard conditional expectation from $\CRG$ to
$\CZGZ$ \scite{\RenaultBook}{II.4.8}, \scite{\Renault}{4.3}.
For $f$ in  $\CCG$ recall that $E(f)$ coincides with the restriction of $f$ to
$\Gz$, so we may write \lcite{\ApplicationOfKeyWithSup} as 
  $$
  \|E(f)\|\leq \|\rep(f)\|
  \for f\in \CCG.
  \subeqmark FinalApplicationOfKey
  $$

Letting $B$ be the range of $\rep$,
we claim that there exists a bounded linear map $F$ from $B$ to
$\CZGZ$ such that the diagram

  \bigskip\beginmypicture
  \setcoordinatesystem units <0.0015truecm, -0.0015truecm> point at 0 0
  \pouquinho = 350
  \put {$\CRG$} at 0000  0000
  \put {$B$} at 2000 0000
  \morph {0200}{0000}{2000}{0000} \put{$\pi$} at 1050 -200
  \put {$\CZGZ$} at 0000  1500
  \morph {0000}{0000}{0000}{1500} \put{$E$} at -250 750
  \pouquinho = 700
  \morph {2000}{0000}{0000}{1500} \put{$F$} at  1350 750
  \endmypicture
  \bigskip
  \noindent commutes.  We first define $F$ on the dense *-subalgebra
$\rep\big(\CCG\big)\c B$, by
  $$
  F\big(\rep(f)\big) = E(f)
  \for f\in \CCG.
  $$
  By \lcite{\FinalApplicationOfKey} this is well defined and bounded,
and hence may be continuously extended to the whole of $B$.  The
extension will then clearly  satisfy the required conditions.

Let $a\in\CRG$ be such that $\rep(a)=0$.  Then 
  $$
  0 = F\big(\rep(a^*)\rep(a)\big) =
  F\big(\rep(a^*a)\big) =
  E(a^*a).
  $$
  Since $E$ is faithful \scite{\Renault}{4.3.ii}, we deduce that
$a=0$, hence concluding the proof of (a).

We now turn to proving (b).  Consider a representation $\rep$ of $\CRG$
whose kernel coincides with $J$.  Such a representation may be
obtained by faithfully embedding $\CRG/J$ as an algebra of operators
on a Hilbert space.

Arguing by contradiction,  if the intersection of $J$ with $\CZGZ$ is
zero, then the restriction  
$\rep|_{C_0(\Gz)}$ is faithful and hence $\rep$ itself is faithful by (a),
from which one would deduce that $J$ is zero.
  \proofend

  \endgroup

  \begingroup
  \bigskip\bigskip 
  \font \sc = cmcsc8 \sc
  \parskip = -1pt

  Departamento de Matem\'atica 

  Universidade Federal de Santa Catarina

  88040-900 -- Florian\'opolis -- Brasil

  \eightrm r@exel.com.br

  \endgroup
  \bye